\pdfoutput=1
\documentclass[12pt]{extarticle}
\usepackage{comment}
\usepackage{caption}
\usepackage{multicol}
\usepackage{xfrac} 
\usepackage{microtype}
\oddsidemargin \evensidemargin
\usepackage{floatrow}
\usepackage[dvipsnames]{xcolor}

\usepackage{amsmath, amssymb, amsfonts}
\usepackage{mathtools}
\usepackage{enumerate}
\usepackage{multicol}
\usepackage[utf8]{inputenc}
\usepackage{bbm}
\usepackage{mathrsfs}
\usepackage{tikz,pgfplots,tikz-3dplot}
\usepackage{changepage}
\usepackage{algorithm}
\usepackage{algpseudocode}
\algdef{SE}[SUBALG]{Indent}{EndIndent}{}{\algorithmicend\ }%
\algtext*{Indent}
\algtext*{EndIndent}
\usepackage{listings}
\usepackage{listings}
\usepackage{xcolor}
\definecolor{maccolor}{rgb}{0.3,0.3,0.8}
\lstdefinelanguage{Macaulay2}
{
basicstyle={\ttfamily},
keywordstyle={\color{maccolor!80!black}},
commentstyle={\color{gray}},
stringstyle={\color{red!40!black}},
rulecolor=\color{maccolor},
basewidth={1.2ex}, 
sensitive=false,
morestring=[b]",
escapechar={`},
escapebegin={\rmfamily},
morekeywords={
ConnectionMatrices, connectionForm, pfaffianSystem, gaugeMatrix, gaugeTransform, holonomicRank, isEpsilonFactorized, isIntegrable, makeWeylAlgebra, standardMonomials, baseFractionField, pfaffianSystem
about,abs,AbstractToricVarieties,accumulate,Acknowledgement,acos,acosh,acot,addCancelTask,addDependencyTask,addEndFunction,addHook,AdditionalPaths,addStartFunction,addStartTask,Adjacent,adjoint,AdjointIdeal,AffineVariety,AfterEval,AfterNoPrint,AfterPrint,agm,AInfinity,alarm,AlgebraicSplines,Algorithm,Alignment,all,AllCodimensions,allowableThreads,ambient,analyticSpread,Analyzer,AnalyzeSheafOnP1,ancestor,ancestors,ANCHOR,and,andP,AngleBarList,ann,annihilator,antipode,any,append,applicationDirectory,applicationDirectorySuffix,apply,applyKeys,applyPairs,applyTable,applyValues,apropos,argument,Array,arXiv,Ascending,ascii,asin,asinh,ass,assert,associatedGradedRing,associatedPrimes,AssociativeAlgebras,AssociativeExpression,atan,atan2,atEndOfFile,Authors,autoload,AuxiliaryFiles,backtrace,Bag,Bareiss,baseFilename,BaseFunction,baseName,baseRing,baseRings,BaseRow,BasicList,basis,BasisElementLimit,Bayer,BeforePrint,beginDocumentation,BeginningMacaulay2,Benchmark,benchmark,Bertini,BesselJ,BesselY,betti,BettiCharacters,BettiTally,between,BGG,BIBasis,Binary,BinaryOperation,Binomial,binomial,BinomialEdgeIdeals,Binomials,BKZ,BlockMatrix,BLOCKQUOTE,BODY,Body,BoijSoederberg,BOLD,Book3264Examples,Boolean,BooleanGB,borel,Boxes,BR,break,Browse,Bruns,cache,CacheExampleOutput,CacheFunction,CacheTable,cacheValue,CallLimit,cancelTask,capture,catch,Caveat,CC,CDATA,ceiling,Center,centerString,Certification,ChainComplex,chainComplex,ChainComplexExtras,ChainComplexMap,ChainComplexOperations,ChangeMatrix,char,CharacteristicClasses,characters,charAnalyzer,check,CheckDocumentation,chi,Chordal,class,Classic,clean,clearAll,clearEcho,clearOutput,close,closeIn,closeOut,ClosestFit,CODE,code,codim,CodimensionLimit,coefficient,CoefficientRing,coefficientRing,coefficients,Cofactor,CohenEngine,CohenTopLevel,CoherentSheaf,CohomCalg,cohomology,coimage,CoincidentRootLoci,coker,cokernel,collectGarbage,columnAdd,columnate,columnMult,columnPermute,columnRankProfile,columnSwap,combine,Command,commandInterpreter,commandLine,COMMENT,commonest,commonRing,comodule,CompactMatrix,compactMatrixForm,CompiledFunction,CompiledFunctionBody,CompiledFunctionClosure,Complement,complement,complete,CompleteIntersection,CompleteIntersectionResolutions,Complexes,ComplexField,components,compose,compositions,compress,concatenate,conductor,ConductorElement,cone,Configuration,ConformalBlocks,conjugate,connectionCount,Consequences,Constant,Constants,constParser,content,continue,contract,Contributors,ConvexInterface,conwayPolynomial,ConwayPolynomials,copy,copyDirectory,copyFile,copyright,Core,CorrespondenceScrolls,cos,cosh,cot,CotangentSchubert,cotangentSheaf,coth,cover,coverMap,cpuTime,createTask,Cremona,csc,csch,current,currentColumnNumber,currentDirectory,currentFileDirectory,currentFileName,currentLayout,currentLineNumber,currentPackage,currentString,currentTime,Cyclotomic,Database,Date,DD,dd,deadParser,debug,debugError,DebuggingMode,debuggingMode,debugLevel,DecomposableSparseSystems,Decompose,decompose,deepSplice,Default,default,defaultPrecision,Degree,degree,degreeLength,DegreeLift,DegreeLimit,DegreeMap,DegreeOrder,DegreeRank,Degrees,degrees,degreesMonoid,degreesRing,delete,demark,denominator,Dense,Density,Depth,depth,Descending,Descent,Describe,describe,Description,det,determinant,DeterminantalRepresentations,DGAlgebras,diagonalMatrix,diameter,Dictionary,dictionary,dictionaryPath,diff,DiffAlg,difference,dim,directSum,disassemble,discriminant,dismiss,Dispatch,distinguished,DIV,Divide,divideByVariable,DivideConquer,DividedPowers,Divisor,DL,Dmodules,do,doc,docExample,docTemplate,document,DocumentTag,Down,drop,DT,dual,eagonNorthcott,EagonResolution,echoOff,echoOn,EdgeIdeals,edit,EigenSolver,eigenvalues,eigenvectors,eint,EisenbudHunekeVasconcelos,elapsedTime,elapsedTiming,elements,Eliminate,eliminate,Elimination,EliminationMatrices,EllipticCurves,EllipticIntegrals,else,EM,Email,End,end,endl,endPackage,Engine,engineDebugLevel,EngineRing,EngineTests,entries,EnumerationCurves,environment,Equation,EquivariantGB,erase,erf,erfc,error,errorDepth,euler,EulerConstant,eulers,even,EXAMPLE,ExampleFiles,ExampleItem,examples,ExampleSystems,Exclude,exec,exit,exp,expectedReesIdeal,expm1,exponents,export,exportFrom,exportMutable,Expression,expression,Ext,extend,ExteriorIdeals,ExteriorModules,exteriorPower,Factor,factor,false,Fano,FastMinors,FastNonminimal,FGLM,File,fileDictionaries,fileExecutable,fileExists,fileExitHooks,fileLength,fileMode,FileName,FilePosition,fileReadable,fileTime,fileWritable,fillMatrix,findFiles,findHeft,FindOne,findProgram,findSynonyms,FiniteFittingIdeals,First,first,firstkey,FirstPackage,fittingIdeal,flagLookup,FlatMonoid,flatten,flattenRing,Flexible,flip,floor,flush,fold,FollowLinks,for,forceGB,fork,FormalGroupLaws,Format,format,formation,FourierMotzkin,FourTiTwo,fpLLL,frac,fraction,FractionField,frames,FrobeniusThresholds,from,fromDividedPowers,fromDual,Function,FunctionApplication,FunctionBody,functionBody,FunctionClosure,FunctionFieldDesingularization,fusePairs,futureParser,GaloisField,Gamma,gb,GBDegrees,gbRemove,gbSnapshot,gbTrace,gcd,gcdCoefficients,gcdLLL,GCstats,genera,GeneralOrderedMonoid,GenerateAssertions,generateAssertions,generator,generators,Generic,GenericInitialIdeal,genericMatrix,genericSkewMatrix,genericSymmetricMatrix,gens,genus,get,getc,getChangeMatrix,getenv,getGlobalSymbol,getNetFile,getNonUnit,getPrimeWithRootOfUnity,getSymbol,getWWW,GF,gfanInterface,Givens,GKMVarieties,GLex,Global,global,globalAssign,globalAssignFunction,GlobalAssignHook,globalAssignment,globalAssignmentHooks,GlobalDictionary,GlobalHookStore,globalReleaseFunction,GlobalReleaseHook,Gorenstein,GradedLieAlgebras,GradedModule,gradedModule,GradedModuleMap,gradedModuleMap,gramm,GraphicalModels,GraphicalModelsMLE,Graphics,graphIdeal,graphRing,Graphs,Grassmannian,GRevLex,GroebnerBasis,groebnerBasis,GroebnerBasisOptions,GroebnerStrata,GroebnerWalk,groupID,GroupLex,GroupRevLex,GTZ,Hadamard,handleInterrupts,HardDegreeLimit,hash,HashTable,hashTable,HEAD,HEADER1,HEADER2,HEADER3,HEADER4,HEADER5,HEADER6,HeaderType,Heading,Headline,Heft,heft,Height,height,help,Hermite,hermite,Hermitian,HH,hh,HigherCIOperators,HighestWeights,Hilbert,hilbertFunction,hilbertPolynomial,hilbertSeries,HodgeIntegrals,hold,Holder,Hom,homeDirectory,HomePage,Homogeneous,Homogeneous2,homogenize,homology,homomorphism,HomotopyLieAlgebra,hooks,horizontalJoin,HorizontalSpace,HR,HREF,HTML,html,httpHeaders,Hybrid,HyperplaneArrangements,Hypertext,hypertext,HypertextContainer,HypertextParagraph,icFracP,icFractions,icMap,icPIdeal,id,Ideal,ideal,idealizer,identity,if,IgnoreExampleErrors,ii,image,imaginaryPart,IMG,ImmutableType,importFrom,in,incomparable,Increment,independentSets,indeterminate,IndeterminateNumber,Index,index,indexComponents,IndexedVariable,IndexedVariableTable,indices,inducedMap,inducesWellDefinedMap,InexactField,InexactFieldFamily,InexactNumber,InfiniteNumber,infinity,info,InfoDirSection,infoHelp,Inhomogeneous,input,Inputs,insert,installAssignmentMethod,installedPackages,installHilbertFunction,installMethod,installMinprimes,installPackage,InstallPrefix,instance,instances,IntegralClosure,integralClosure,integrate,IntermediateMarkUpType,interpreterDepth,intersect,intersectInP,Intersection,intersection,interval,InvariantRing,inverse,InverseMethod,inversePermutation,Inverses,inverseSystem,InverseSystems,Invertible,InvolutiveBases,irreducibleCharacteristicSeries,irreducibleDecomposition,isAffineRing,isANumber,isBorel,isCanceled,isCommutative,isConstant,isDirectory,isDirectSum,isEmpty,isField,isFinite,isFinitePrimeField,isFreeModule,isGlobalSymbol,isHomogeneous,isIdeal,isInfinite,isInjective,isInputFile,isIsomorphism,isLinearType,isListener,isLLL,isMember,isModule,isMonomialIdeal,isNormal,isOpen,isOutputFile,isPolynomialRing,isPrimary,isPrime,isPrimitive,isPseudoprime,isQuotientModule,isQuotientOf,isQuotientRing,isReady,isReal,isReduction,isRegularFile,isRing,isSkewCommutative,isSorted,isSquareFree,isStandardGradedPolynomialRing,isSubmodule,isSubquotient,isSubset,isSupportedInZeroLocus,isSurjective,isTable,isUnit,isWellDefined,isWeylAlgebra,ITALIC,Iterate,Jacobian,jacobian,jacobianDual,Jets,Join,join,Jupyter,K3Carpets,K3Surfaces,Keep,KeepFiles,KeepZeroes,ker,kernel,kernelLLL,kernelOfLocalization,Key,keys,Keyword,Keywords,kill,koszul,Kronecker,KustinMiller,LABEL,last,lastMatch,LATER,LatticePolytopes,Layout,lcm,leadCoefficient,leadComponent,leadMonomial,leadTerm,Left,left,length,LengthLimit,letterParser,Lex,LexIdeals,LI,Licenses,LieTypes,lift,liftable,Limit,limitFiles,limitProcesses,Linear,LinearAlgebra,LinearTruncations,lineNumber,lines,LINK,linkFile,List,list,listForm,listLocalSymbols,listSymbols,listUserSymbols,LITERAL,LLL,LLLBases,lngamma,load,loadDepth,LoadDocumentation,loadedFiles,loadedPackages,loadPackage,Local,local,localDictionaries,LocalDictionary,localize,LocalRings,locate,log,log1p,LongPolynomial,lookup,lookupCount,LowerBound,LUdecomposition,M0nbar,M2CODE,Macaulay2Doc,makeDirectory,MakeDocumentation,makeDocumentTag,MakeHTML,MakeInfo,MakeLinks,makePackageIndex,MakePDF,makeS2,Manipulator,map,MapExpression,MapleInterface,markedGB,Markov,MarkUpType,match,mathML,Matrix,matrix,MatrixExpression,Matroids,max,maxAllowableThreads,maxExponent,MaximalRank,maxPosition,MaxReductionCount,MCMApproximations,member,memoize,memoizeClear,memoizeValues,MENU,merge,mergePairs,META,method,MethodFunction,MethodFunctionBinary,MethodFunctionSingle,MethodFunctionWithOptions,methodOptions,methods,midpoint,min,minExponent,mingens,mingle,minimalBetti,MinimalGenerators,MinimalMatrix,minimalPresentation,minimalPresentationMap,minimalPresentationMapInv,MinimalPrimes,minimalPrimes,minimalReduction,Minimize,minimizeFilename,MinimumVersion,minors,minPosition,minPres,minprimes,Minus,minus,Miura,MixedMultiplicity,mkdir,mod,Module,module,ModuleDeformations,modulo,MonodromySolver,Monoid,monoid,MonoidElement,Monomial,MonomialAlgebras,monomialCurveIdeal,MonomialIdeal,monomialIdeal,MonomialIntegerPrograms,MonomialOrbits,MonomialOrder,Monomials,monomials,MonomialSize,monomialSubideal,moveFile,multidegree,multidoc,multigraded,MultigradedBettiTally,MultiGradedRationalMap,multiplicity,MultiplicitySequence,MultiplierIdeals,MultiplierIdealsDim2,MultiprojectiveVarieties,mutable,MutableHashTable,mutableIdentity,MutableList,MutableMatrix,mutableMatrix,NAGtypes,Name,nanosleep,Nauty,NautyGraphs,NCAlgebra,NCLex,needs,needsPackage,Net,net,NetFile,netList,new,newClass,newCoordinateSystem,NewFromMethod,newline,NewMethod,newNetFile,NewOfFromMethod,NewOfMethod,newPackage,newRing,nextkey,nextPrime,nil,NNParser,NoetherianOperators,NoetherNormalization,NonminimalComplexes,nonspaceAnalyzer,NoPrint,norm,normalCone,Normaliz,NormalToricVarieties,not,Nothing,notify,notImplemented,NTL,null,nullaryMethods,nullhomotopy,nullParser,nullSpace,Number,number,NumberedVerticalList,numcols,numColumns,numerator,numeric,NumericalAlgebraicGeometry,NumericalCertification,NumericalImplicitization,NumericalLinearAlgebra,NumericalSchubertCalculus,numericInterval,NumericSolutions,numgens,numRows,numrows,odd,oeis,of,ofClass,OL,OldPolyhedra,OldToricVectorBundles,on,OneExpression,OnlineLookup,OO,oo,ooo,oooo,openDatabase,openDatabaseOut,openFiles,openIn,openInOut,openListener,OpenMath,openOut,openOutAppend,operatorAttributes,Option,OptionalComponentsPresent,optionalSignParser,Options,options,OptionTable,optP,or,Order,order,OrderedMonoid,orP,OutputDictionary,Outputs,override,pack,Package,package,PackageCitations,PackageDictionary,PackageExports,PackageImports,PackageTemplate,packageTemplate,pad,pager,PairLimit,pairs,PairsRemaining,PARA,Parametrization,parent,Parenthesize,Parser,Parsing,part,Partition,partition,partitions,parts,path,pdim,peek,PencilsOfQuadrics,Permanents,permanents,permutations,pfaffians,PHCpack,PhylogeneticTrees,pi,PieriMaps,pivots,PlaneCurveSingularities,plus,poincare,poincareN,Points,polarize,poly,Polyhedra,Polymake,PolynomialRing,Posets,Position,position,positions,PositivityToricBundles,POSIX,Postfix,Power,power,powermod,PRE,Precision,precision,Prefix,prefixDirectory,prefixPath,preimage,prepend,presentation,pretty,primaryComponent,PrimaryDecomposition,primaryDecomposition,PrimaryTag,PrimitiveElement,Print,print,printerr,printingAccuracy,printingLeadLimit,printingPrecision,printingSeparator,printingTimeLimit,printingTrailLimit,printString,printWidth,processID,Product,product,ProductOrder,profile,profileSummary,Program,programPaths,ProgramRun,Proj,Projective,ProjectiveHilbertPolynomial,projectiveHilbertPolynomial,ProjectiveVariety,promote,protect,Prune,prune,PruneComplex,pruningMap,Pseudocode,pseudocode,pseudoRemainder,Pullback,PushForward,pushForward,Python,QQ,QQParser,QRDecomposition,QthPower,Quasidegrees,QuaternaryQuartics,QuillenSuslin,quit,Quotient,quotient,quotientRemainder,QuotientRing,Radical,radical,RadicalCodim1,radicalContainment,RaiseError,random,RandomCanonicalCurves,RandomComplexes,RandomCurves,RandomCurvesOverVerySmallFiniteFields,RandomGenus14Curves,RandomIdeals,randomKRationalPoint,RandomMonomialIdeals,randomMutableMatrix,RandomObjects,RandomPlaneCurves,RandomPoints,RandomSpaceCurves,Range,rank,RationalMaps,RationalPoints,RationalPoints2,ReactionNetworks,read,readDirectory,readlink,readPackage,RealField,RealFP,realPart,realpath,RealQP,RealQP1,RealRoots,RealRR,RealXD,recursionDepth,recursionLimit,Reduce,reducedRowEchelonForm,reduceHilbert,reductionNumber,ReesAlgebra,reesAlgebra,reesAlgebraIdeal,reesIdeal,References,ReflexivePolytopesDB,regex,regexQuote,registerFinalizer,regSeqInIdeal,Regularity,regularity,relations,RelativeCanonicalResolution,relativizeFilename,Reload,remainder,RemakeAllDocumentation,remove,removeDirectory,removeFile,removeLowestDimension,reorganize,replace,RerunExamples,res,reshape,ResidualIntersections,ResLengthThree,Resolution,resolution,ResolutionsOfStanleyReisnerRings,restart,Result,resultant,Resultants,return,returnCode,Reverse,reverse,RevLex,Right,right,Ring,ring,RingElement,RingFamily,ringFromFractions,RingMap,rootPath,roots,rootURI,rotate,round,rowAdd,RowExpression,rowMult,rowPermute,rowRankProfile,rowSwap,RR,RRi,rsort,run,RunDirectory,RunExamples,RunExternalM2,runHooks,runLengthEncode,runProgram,same,saturate,Saturation,scan,scanKeys,scanLines,scanPairs,scanValues,schedule,schreyerOrder,Schubert,Schubert2,SchurComplexes,SchurFunctors,SchurRings,SCRIPT,scriptCommandLine,ScriptedFunctor,SCSCP,searchPath,sec,sech,SectionRing,SeeAlso,seeParsing,SegreClasses,select,selectInSubring,selectVariables,SelfInitializingType,SemidefiniteProgramming,Seminormalization,separate,SeparateExec,separateRegexp,Sequence,sequence,Serialization,serialNumber,Set,set,setEcho,setGroupID,setIOExclusive,setIOSynchronized,setIOUnSynchronized,setRandomSeed,setup,setupEmacs,sheaf,SheafExpression,sheafExt,sheafHom,SheafOfRings,shield,ShimoyamaYokoyama,short,show,showClassStructure,showHtml,showStructure,showTex,showUserStructure,SimpleDoc,simpleDocFrob,SimplicialComplexes,SimplicialDecomposability,SimplicialPosets,SimplifyFractions,sin,singularLocus,sinh,size,size2,SizeLimit,SkewCommutative,SlackIdeals,sleep,SLnEquivariantMatrices,SLPexpressions,SMALL,smithNormalForm,solve,someTerms,Sort,sort,sortColumns,SortStrategy,source,SourceCode,SourceRing,SPACE,SpaceCurves,SPAN,span,SparseMonomialVectorExpression,SparseResultants,SparseVectorExpression,Spec,SpechtModule,SpecialFanoFourfolds,specialFiber,specialFiberIdeal,SpectralSequences,splice,splitWWW,sqrt,SRdeformations,stack,stacksProject,Standard,standardForm,standardPairs,StartWithOneMinor,stashValue,StatePolytope,StatGraphs,status,stderr,stdio,step,StopBeforeComputation,stopIfError,StopWithMinimalGenerators,Strategy,String,STRONG,StronglyStableIdeals,STYLE,Style,style,SUB,sub,SubalgebraBases,sublists,submatrix,submatrixByDegrees,Subnodes,subquotient,SubringLimit,Subscript,subscript,SUBSECTION,subsets,substitute,substring,subtable,Sugarless,Sum,sum,SumOfTwists,SumsOfSquares,SUP,super,SuperLinearAlgebra,Superscript,superscript,support,SVD,SVDComplexes,switch,SwitchingFields,sylvesterMatrix,Symbol,symbol,SymbolBody,symbolBody,SymbolicPowers,symlinkDirectory,symlinkFile,symmetricAlgebra,symmetricAlgebraIdeal,symmetricKernel,SymmetricPolynomials,symmetricPower,synonym,SYNOPSIS,syz,Syzygies,SyzygyLimit,SyzygyMatrix,SyzygyRows,syzygyScheme,TABLE,Table,table,take,Tally,tally,tan,TangentCone,tangentCone,tangentSheaf,tanh,target,Task,taskResult,TateOnProducts,TD,temporaryFileName,tensor,tensorAssociativity,TensorComplexes,terminalParser,terms,TEST,Test,testExample,testHunekeQuestion,TestIdeals,TestInput,tests,TEX,tex,TeXmacs,texMath,Text,TH,then,Thing,ThinSincereQuivers,ThreadedGB,threadVariable,Threshold,throw,Time,time,times,timing,TITLE,TO,to,TO2,toAbsolutePath,toCC,toDividedPowers,toDual,toExternalString,toField,TOH,toList,toLower,top,top,topCoefficients,Topcom,topComponents,topLevelMode,Tor,TorAlgebra,Toric,ToricInvariants,ToricTopology,ToricVectorBundles,toRR,toRRi,toSequence,toString,TotalPairs,toUpper,TR,trace,transpose,TriangularSets,Tries,Trim,trim,Triplets,Tropical,true,Truncate,truncate,truncateOutput,Truncations,try,TSpreadIdeals,TT,tutorial,Type,TypicalValue,typicalValues,UL,ultimate,unbag,uncurry,Undo,undocumented,uniform,uninstallAllPackages,uninstallPackage,Unique,unique,Units,Unmixed,unsequence,unstack,Up,UpdateOnly,UpperTriangular,URL,urlEncode,Usage,use,UseCachedExampleOutput,UseHilbertFunction,UserMode,userSymbols,UseSyzygies,utf8,utf8check,validate,value,values,Variable,VariableBaseName,Variables,Variety,variety,vars,Vasconcelos,Vector,vector,VectorExpression,VectorFields,VectorGraphics,Verbose,Verbosity,Verify,VersalDeformations,versalEmbedding,Version,version,VerticalList,VerticalSpace,viewHelp,VirtualResolutions,VirtualTally,VisibleList,Visualize,wait,WebApp,wedgeProduct,weightRange,Weights,WeylAlgebra,WeylGroups,when,whichGm,while,width,wikipedia,Wrap,wrap,WrapperType,XML,xor,youngest,zero,ZeroExpression,zeta,ZZ,ZZParser}
}
\lstalias{Macaulay2output}{Macaulay2}

\usepackage{titlesec}
\titleformat{\subsubsection}[runin]{\normalfont\bfseries}{\thesubsubsection}{0.5em}{}[]

\renewenvironment{thebibliography}[1]{
  \begin{oldthebibliography}{#1}
    \setlength{\itemsep}{0.5em}
    \setlength{\parskip}{0em}
}
{
  \end{oldthebibliography}
}

\usepackage[english]{babel}

\usepackage[colorlinks=true, allcolors=blue,linkcolor=MidnightBlue, citecolor=Blue,backref=page,urlcolor=black]{hyperref}

\usepackage{amsthm}
\usepackage{cleveref}

\newtheorem{theorem}{Theorem}[section]

\newtheorem{proposition}[theorem]{Proposition}

\theoremstyle{definition}

\newtheorem{definition}[theorem]{Definition}

\newtheorem{examplex}[theorem]{Example}
\newenvironment{example}
{\pushQED{\qed}\examplex}
{\popQED\endexamplex}

\newtheorem{remarkx}[theorem]{Remark}
\newenvironment{remark}
{\pushQED{\qed}\remarkx}
{\popQED\endremarkx}

\numberwithin{equation}{section}

\newtheoremstyle{citing}
{}
{}
{\itshape}
{}
{\bfseries}
{\textbf{.}}
{.5em}
{\thmnote{#3}}
{\theoremstyle{citing}
}

\DeclareMathOperator{\lt}{lt}

\DeclareMathOperator{\NN}{\mathbb{N}}
\DeclareMathOperator{\CC}{\mathbb{C}}
\DeclareMathOperator{\RR}{\mathbb{R}}

\DeclareMathOperator{\init}{in}
\DeclareMathOperator{\leadTerm}{lt}

\DeclareMathOperator{\rank}{rank}

\DeclareMathOperator{\GL}{GL}

\DeclareMathOperator{\Sol}{Sol}

\renewcommand{\d}{\mathrm{d}}

\pgfplotsset{compat=1.18}

\title{Connection Matrices in {\em Macaulay2}}
\author{Paul Görlach${}^\sharp$, Joris Koefler${}^\flat$, Anna-Laura Sattelberger${}^\flat$,
 Mahrud Sayrafi${}^{\circ}$, \\Hendrik Schroeder${}^\diamond$,
Nicolas Weiss${}^\flat$, and Francesca Zaffalon${}^{\flat\vartriangle}$}
\date{}

\begin{document}
\maketitle
\thispagestyle{empty}

\begin{abstract}
In this article, we describe the theoretical foundations of the {\em Macaulay2} package {\tt ConnectionMatrices} and explain how to use it. For a left ideal in the Weyl algebra that is of finite holonomic rank, we implement the computation of the encoded system of linear PDEs in connection form with respect to an elimination term order that depends on a chosen positive weight vector. We also implement the gauge transformation for carrying out a change of basis over the field of rational functions. We demonstrate all implemented algorithms with examples.
\end{abstract}

{\hypersetup{linkcolor=black}
\setcounter{tocdepth}{2}
\renewcommand{\baselinestretch}{.9}\normalsize
\tableofcontents
\renewcommand{\baselinestretch}{1.0}\normalsize
}

\vfill

{\small
\noindent ${}^\sharp$ Otto von Guericke University Magdeburg, Germany \hfill  {\tt paul.goerlach@ovgu.de}

\noindent ${}^{\flat}$  Max Planck Institute for Mathematics in the Sciences, Leipzig, Germany  \hfill {\tt $\{$ joris.koefler,} \\ \hspace*{1.4mm} {\tt anna-laura.sattelberger, nicolas.weiss,} and {\tt francesca.zaffalon $\}$ @ mis.mpg.de}

\noindent ${}^\circ$ Fields Institute, Toronto, Canada \hfill {\tt mahrud@fields.utoronto.ca}

\noindent ${}^\diamond$ Technical University of Berlin, Germany \hfill {\tt h.schroeder@tu-berlin.de}

\noindent \hspace*{-.5mm}${}^\vartriangle$ Weizmann Institute of Science, Israel}

\newpage
\addcontentsline{toc}{section}{Introduction}
\section*{Introduction}
Systems of homogeneous, linear partial differential equations (PDEs) with polynomial coefficients are encoded by left ideals in the Weyl algebra, denoted $D_n=\CC[x_1,\ldots,x_n]\langle \partial_1,\ldots,\partial_n \rangle$ (or just~$D$).
Such systems can be systematically written as a first-order matrix system, i.e., in ``connection form'' $\d - A\wedge \, $, by utilizing Gröbner bases in the Weyl algebra~\cite{SST00}. In the case of a single ordinary differential equation (ODE), 
the respective matrix is known under the name of ``companion matrix.'' {
Consider for instance the $D_1$-ideal $I = \, \left\langle x\partial^3-(x+1)\partial + 1 \right\rangle$. It represents the third-order ODE $\,x \cdot f'''(x)-(x+1) \cdot f'(x)+1 = 0\,$ whose companion matrix~is
\begin{align*}
P \,=\, \begin{pmatrix}
            0 & 1 & 0\\
            0 & 0 & 1\\
            -\frac{1}{x} & \frac{x+1}{x} & 0
        \end{pmatrix}.
\end{align*}
Then for any solution $f$ to the considered ODE and any non-zero 
number $u$, one has 
\begin{align*}
    \begin{pmatrix}
        f'(u) \\
        f''(u)\\
        f'''(u)
    \end{pmatrix} \,=\,
        P(u) \cdot 
            \begin{pmatrix}
                f(u)\\
                f'(u)\\
                f''(u)
            \end{pmatrix} .
\end{align*}
Connection matrices generalize the concept of companion matrices of ODEs to homogeneous, linear PDEs with polynomial coefficients. They are an important tool for numerical algorithms, such as the holonomic gradient descent~\cite{HGD} for the optimization of real-valued holonomic functions.}
The systematic computation of connection matrices in software requires Gröbner bases in the rational Weyl algebra $R_n=\CC(x_1,\ldots,x_n)\langle \partial_1,\ldots,\partial_n \rangle$. These, however, are not yet available in the $D$-module packages in standard open-source computer algebra software.

In this article, we provide the theoretical background for the {\em Macaulay2}~\cite{M2} package {\tt ConnectionMatrices}. The implemented functionalities include the computation of connection matrices of $D$-ideals for elimination term orders on the Weyl algebra with respect to positive weight vectors. In particular, this required the implementation of the normal form algorithm over the rational Weyl algebra. We also implemented the gauge transformation for carrying out changes of basis over the field of rational functions. While the theory presented here is formulated over the complex numbers as the field of coefficients, the implementations are---as is usual---over the rational numbers. We also allow for the dependence on parameters, such as a ``small parameter''~$\varepsilon$, as is commonly used for the dimensional regularization of Feynman integrals in particle physics.

\smallskip
For our implementations, we make use of the {\em Macaulay2} package {\tt WeylAlgebras}~\cite{WeylAlgebra}, which is included in the package collection~{\tt Dmodules}~\cite{DmodM2}. Our package was first made available via the MathRepo~\cite{mathrepo} hosted by MPI~MiS at \url{https://mathrepo.mis.mpg.de/ConnectionMatrices}.
It follows the FAIR data principles of the mathematical research data initiative \href{https://www.mardi4nfdi.de/about/mission}{MaRDI}~\cite{Mardi}, which aim to improve {\bf f}indability, {\bf a}ccessibility, {\bf i}nteroperability, and {\bf r}euse of digital~assets. By now, \href{https://macaulay2.com/doc/Macaulay2-1.25.05/share/doc/Macaulay2/ConnectionMatrices/html/index.html}{\tt ConnectionMatrices} is contained as a package in the {\em Macaulay2} release 1.25.05, together with a documentation of all implemented commands.

\pagebreak
\noindent{\bf Notation.} 
Elements of the (rational) Weyl algebra are typically denoted by the letter $P \in D_n$ (or $P\in R_n$). Left ideals in the Weyl algebra are denoted by $I=\langle P_1,\ldots,P_k\rangle\subset D_n$, $m=\rank(I)$ denotes their holonomic rank, and $G=\{ G_1,\ldots,G_{\ell} \}$ denotes Gröbner bases. The letters $A_i$ denote the connection matrices of a $D$-ideal~$I$; they are $m\times m$ matrices with entries in the field of rational functions. Equivalently, one can encode the $A_i$'s in a single $m\times m$ matrix~$A$ of differential one-forms, the connection matrix of~$I$. We use $\{s_1,\ldots,s_m\}$ to denote a $\CC(x_1,\ldots,x_n)$-basis of $R_n/R_nI$. In our implementations, the $s_i$'s are typically chosen as the standard monomials of a Gröbner basis of $R_nI$. A gauge transformation of the connection matrices is encoded by an invertible matrix $g\in \operatorname{GL}_m(\CC(x_1,\ldots,x_n))$.

\bigskip
\noindent{\bf Outline.} \Cref{sec:WeylAlgebra} recalls background on Gröbner bases in non-commutative rings of differential operators. \Cref{sec:connection} explains how to systematically write systems of linear PDEs in connection form. In particular, it presents the implemented algorithms in pseudo-code. \Cref{sec:implementations} explains the implemented functionalities of the package {\tt ConnectionMatrices} and demonstrates all commands via examples.

\section{Gröbner bases in Weyl algebras}\label{sec:WeylAlgebra}
\subsection{The (rational) Weyl algebra}
Homogeneous, linear partial differential equations with polynomial coefficients are encoded as linear differential operators. These are elements of the ($n$-th) Weyl algebra $D_n$ (or just~$D$),
\begin{align}
 D_n \, \coloneqq \, \CC[x_1,\ldots,x_n] \langle \partial_1,\ldots,\partial_n \rangle  ,
\end{align}
which denotes the non-commutative $\CC$-algebra obtained from the free $\CC$-algebra generated by $x_1,\ldots,x_n$ and $\partial_1,\ldots,\partial_n$ by imposing the following relations: all generators are assumed to commute, except $x_i$ and $\partial_i$. Their commutator obeys Leibniz' rule: $[\partial_i,x_i]=1$, $i=1,\ldots,n$. Each element $P\in D_n$ can be uniquely expressed as
\begin{align}\label{eq:normord}
    P \,=\, \sum_{(\alpha,\beta) \,\in \, E} c_{\alpha,\beta}x^{\alpha}\partial^{\beta},
\end{align}
where $E\subset \NN^{2n}$ is a finite set, $c_{\alpha,\beta}\in \CC\setminus \{0\}$, and multi-index notation is used.
We will denote the action of a differential operator on a function $f(x_1,\ldots,x_n)$ by a bullet; for instance, $\partial_i\bullet f = \frac{\partial f}{\partial x_i}$, so that the PDE associated to~\eqref{eq:normord}, 
\begin{align}
\sum_{(\alpha,\beta) \,\in \, E} c_{\alpha,\beta}x_1^{\alpha _1}\cdots x_n^{\alpha_n}f^{(\beta_1,\ldots,\beta_n)}(x_1,\ldots,x_n) \, = \, 0 \, \quad \text{for all } x \, ,
\end{align}
reads $P\bullet f=0$.
A system of PDEs of the form
$\{ P_1\bullet f =0, \, P_2\bullet f=0, \cdots, P_k \bullet f =0\}$
is encoded by the left $D_n$-ideal generated by $P_1,\ldots,P_k$, which we denote by $\langle P_1,\ldots,P_k \rangle\subset D_n$.
Linear differential operators with rational functions as coefficients are elements of the ($n$-th) {\em rational Weyl algebra}, which is denoted by
\begin{align}
 R_n \, \coloneqq \, \CC(x_1,\ldots,x_n) \langle \partial_1,\ldots,\partial_n \rangle \, ,
\end{align}
with the corresponding commutator relations. If clear from the context, we sometimes denote $\CC(x)=\CC(x_1,\ldots,x_n)$ for brevity. When considering ideals in the (rational) Weyl algebra, we always mean left ideals. For an ideal $I\subset D_n$, the {\em holonomic rank} of $I$ is the dimension of the underlying $\CC(x)$-vector space of $R_n/R_nI$. In symbols,
\begin{align}\label{eq:rankI}
    \rank(I) \,\coloneqq\, \dim_{\CC(x)} \left( R_n/R_nI \right) .
\end{align}
Our definition differs from the one given in \cite[Definition~1.4.8]{SST00}; we discuss their equivalence in \Cref{rem:holrank}.
On simply connected domains in $\CC^n$ avoiding the singular locus of $I$, $\rank(I)$ is the dimension of the $\CC$-vector space of holomorphic solutions to $I$, which is 
denoted $\Sol(I)$. This statement follows from the theorem of Cauchy--Kovalevskaya--Kashiwara. 
{The number in~\eqref{eq:rankI} is guaranteed to be finite if $I$ (or its associated $D$-module $D/I$) is {\em holonomic}, see e.g.~\cite[Definition 2.3.6]{HTT08}.
Holonomicity is a central notion for $D$-ideals and it is stronger than having finite holonomic rank.
The two concepts agree for {\em Weyl-closed} $D$-ideals, i.e., $D_n$-ideals $I$ such that $I = R_n I \cap D_n$,
see \cite[Theorem 1.4.15]{SST00} for a more detailed discussion.}

\smallskip

The differences between a $D_n$-ideal $I$ and the $R_n$-ideal $R_nI$ are subtle. However, the connection form of a $D_n$-ideal cannot distinguish between $D_nI$ and $R_nI$. In fact, it depends only on the choice of a $\CC(x)$-basis of $R_n/R_nI$. In the next subsection, we discuss term orders and
Gröbner bases in $D_n$ and $R_n$, mainly following the book~\cite{SST00}.

\subsection{Term orders}\label{sec:Groebner}
We will need to consider total orders on the set of  monomials $\{x^\alpha \partial^\beta\}$ in the Weyl algebra.
Such an order is a {\em  multiplicative monomial order} if
\begin{enumerate}[(1)]
    \item $1\prec x_i\partial_i$ for all $i=1,\ldots,n$ and
    \item $x^\alpha \partial^\beta \prec x^a \partial^b$ implies $x^{\alpha +s}\partial^{\beta+t} \prec x^{a+s} \partial^{b+t}$ for all $(s,t)\in \NN^{2n}$. 
\end{enumerate}
For a fixed multiplicative monomial order $\prec$ on $D_n$ the {\em initial monomial} of $P\in D_n$, denoted $\init_{\prec}(P)$, is the monomial $x^{\alpha}\xi^\beta$ in the (commutative) polynomial ring $ \CC[x_1,\ldots,x_n,\xi_1,\ldots,\xi_n]$ for which $x^\alpha\partial^\beta$ in \eqref{eq:normord} is the largest monomial. The first condition above ensures the compatibility $\init_{\prec}(PQ)=\init_{\prec}(P)\cdot \init_{\prec}(Q)$ with multiplication. The {\em initial ideal} of $I\subset D_n$ with respect to  $\prec$ is the monomial ideal in $\CC[x,\xi]$ that is generated by $\{ \init_{\prec}(P) \mid P\in D_n\}$. A finite set $G = \{G_1,\dots,G_\ell\}\subset D_n$ is a {\em Gröbner basis} of $I$ with respect to $\prec$ if $I=D_nG$ and $\init_{\prec}(I)$ is generated by $\{ \init_{\prec}(G_i) \mid G_i \in G\}$. The \emph{standard monomials} of $I$ with respect to $\prec$ is the set of monomials $x^\alpha \partial^\beta$ which are not contained in $\init_{\prec}(I)$.

\smallskip
A multiplicative monomial order is a {\em term order} on $D_n$ if $1=x^0\partial^0$ is the smallest element of~$\prec$. 
Henceforth, we focus on term orders that arise as the refinement of a partial order given by a weight vector. Weight vectors for $D_n$ are allowed to be taken from the set
\begin{align}\label{eq:weights}
\mathcal{W} \,\coloneqq\,    \left\{(u,v)\in \RR^n\times \RR^n \mid u_i+v_i\geq 0 \ \ \text{for all } i=1,\ldots,n \right\} \, \subset \, \RR^{2n},
\end{align}
i.e., one assigns weight $u_i$ to $x_i$ and weight $v_i$ to $\partial_i$. Each such weight vector induces an increasing, exhaustive filtration of~$D_n$ via the $(u,v)$-weight of differential operators. Later on, we will focus on weight vectors of the form $w=(0,v)$ with $v\in \RR_{>0}^n$ strictly~positive. For $u=0\in \RR^n$ and $v\in \RR^n$ the all-one vector, we denote the resulting weight as $(0,1)\in \RR^{2n}$.

\begin{definition}\label{def:weightref}
 Let $w=(u,v)\in \mathcal{W}$ and $\prec$ be any term order on $D_n$. The order $\prec_{(u,v)}$ is the multiplicative monomial order defined as follows:
\begin{align*}\begin{split}
x^\alpha \partial^\beta  \prec_{(u,v)} x^a \partial^b \, \Leftrightarrow \, \alpha u+\beta v < au+bv \, \text{ or } \, \left(\alpha u+\beta v = au+bv \, \text{ and } \, x^\alpha\partial^\beta \prec x^a\partial^b\right) \, .
\end{split}\end{align*}
\end{definition} 
\noindent That is to say, $\prec$ is used as a tiebreaker in case two monomials have the same $(u,v)$-weight.
This defines a term order if and only if $(u,v)$ is non-negative.

\begin{definition}
A term order $\prec$ on $D_n$ is an {\em elimination term order} if $\partial^\beta \prec \partial^\gamma$ implies $x^\alpha \partial^\beta \prec \partial^\gamma$ for all $\alpha \in \NN^n$.
\end{definition}

Analogous definitions can be made for the rational Weyl algebra $R_n$ with adequate changes to treat $x_1,\dots,x_n$ not as variables, but rather as coefficients: A multiplicative monomial order $\prec'$ on~$R_n$ is a total order on the set of monomials $\{\partial^\beta\}$ with $\partial^\beta \prec' \partial^b$ implying $\partial^{\beta+t} \prec' \partial^{b+t}$, and is called a \emph{term order} if $1 = \partial^0$ is minimal. For $P = \sum_{\beta \in E'} c_\beta(x) \partial^\beta \in R_n$ (with $E' \subset \mathbb{N}^n$, $c_\beta \in \mathbb{C}(x)\setminus \{0\}$), its initial monomial $\init_{\prec'}(P)$ is the element $\xi^\beta \in \mathbb{C}(x)[\xi]$ corresponding to the largest monomial w.r.t.\ $\prec'$ in $\{\partial^\beta \mid \beta \in E'\}$. For an $R_n$-ideal $J$, a finite set $G' = \{G'_1,\dots,G'_\ell\} \subseteq R_n$ is a \emph{Gröbner basis} of $J$ w.r.t.\ $\prec'$ if $J = R_n G'$ and if the initial ideal $\init_{\prec'}(J) \coloneqq  \langle \init_{\prec'}(P) \mid P \in J \rangle \subset \mathbb{C}(x)[\xi]$ is generated by $\{\init_{\prec'}(G'_i) \mid G'_i \in G'\}$. The set of \emph{standard monomials} of $J$ w.r.t.\ $\prec'$ is the set of monomials $\partial^\beta$ with $\xi^\beta \notin \init_{\prec'}(J)$. 

The standard monomials of $R_nI$ with respect to a term order on $R_n$ are a $\CC(x)$-basis of~$R_n/R_nI$---indeed, linear independence over $\CC(x)$ follows directly from the definition, and the normal form algorithm (Algorithm~\ref{alg:normal form algorithm} below) expresses any element of $R_n/R_nI$ as a \mbox{$\CC(x)$-linear} combination of the standard monomials.
Hence, the holonomic rank of a \mbox{$D_n$-ideal} $I$ is the number of standard monomials of~$R_nI$ and may equivalently be expressed as
\begin{align}\label{eq:holrank}
    \rank(I) \,=\, \dim_{\CC(x)} \left( \CC(x)[\xi] /  \init_{\prec'}\left(R_nI\right) \right),
\end{align}
see also \cite[Lemma 1.4.11]{SST00}. Note that $\init_{\prec'}(R_nI) =\CC(x)[\xi]\init_{\prec}(I)$ for $\prec$ an elimination term order on $D_n$ that restricts to $\prec'$. 

\begin{remark}\label{rem:holrank}
  In~\cite{SST00}, the holonomic rank of a $D_n$-ideal $I$ is defined as the dimension of the $\CC(x)$-vector space $\CC(x)[\xi] / \CC(x)[\xi]\cdot \init_{(0,1)}(I)$, where $\init_{(0,1)}(I)$ is the initial ideal of $I$ with respect to the weight vector $(0,1)$. This turns out to be the same as for $\prec_{(0,1)}$, with $\prec$ an arbitrary term order on~$D_n$, see~\cite[Theorem 1.1.6]{SST00}, and hence is a special case of~\eqref{eq:holrank}. 
\end{remark}

In order to construct bases of~$R_n/R_nI$ systematically, we thus compute Gröbner bases in the rational Weyl algebra. 
Let $\prec$ be a term order on $D_n$. We will denote by $\prec'$ its restriction to monomials in the~$\partial_i$'s; this is a term order on~$R_n$.
For any choice of elimination term order $\prec$ on~$D_n$, all $\prec_{(0,v)}$ with strictly positive $v\in \RR_{>0}^n$ are elimination term orders on~$D_n$. 
Stating more clearly than in the last paragraph of~\cite[p.~33]{SST00}, the refinement of the $(0,1)$-weight with respect to an arbitrary term order on $D_n$ does, in general, not result in an elimination term order. For our implementations, we focus on elimination term orders of the form $\prec_{(0,v)}$, with $\prec$ being the lexicographic term order built upon $\partial_1 \succ \cdots \succ \partial_n \succ x_1 \succ \cdots \succ x_n$.

For an $R_n$-ideal $J$, eliminating denominators of a generating set of $J$ leads to the presentation of $J$ as $J = R_n I$ with $I$ being a $D_n$-ideal. We use this small workaround to compute Gröbner bases as follows. 

\begin{proposition}[{\cite[Proposition 1.4.13]{SST00}}]\label{prop1413}
If~$G$ is a Gröbner basis of a $D_n$-ideal~$I$ with respect to an elimination term order $\prec$ on~$D_n$, then $G$ is also a Gröbner basis of the $R_n$-ideal  $R_nI$ with respect to the order~$\prec'$. 
\end{proposition}

\begin{remark}
One could compute the Gröbner bases directly in~$R_n$. However, we refrain from this due to the expected computational overhead caused by 
the bookkeeping and differentiation of rational-function-coefficients, which are necessary to compute S-pairs. Already in the commutative case, we experimentally found Gröbner basis computations to be significantly slower when taking place in $\mathbb{Q}(x_1,\ldots,x_n)[y_1,\ldots, y_n]$ rather than in $\mathbb{Q}[x_1,\ldots, x_n, y_1,\ldots, y_n]$. It is conceivable that there exist specific classes of examples that would benefit from an implementation of Gröbner bases in $R_n$ directly; however, it is unclear how frequently they occur in practical contexts. Our expectation is that our approach of basing the implementations on the already established and optimized Gröbner bases methods in the Weyl algebra to the largest possible extent, turns out favorably in most~circumstances.
\end{remark}

\section{Connection matrices}\label{sec:connection}
Throughout this section, $I$ denotes a $D_n$-ideal of  holonomic rank~$m < \infty$. We explain how to write the system of PDEs encoded by $I$ in matrix form. This can be done in terms of a single $m\times m$  matrix $A$ of differential one-forms, called the {\em connection matrix}, or dually, by $n$ matrices $A_1,\ldots,A_n$ with entries in $\CC(x_1,\ldots,x_n)$.
They depend on the choice of a $\CC(x_1,\ldots,x_n)$-basis of~$R_n/R_nI$, which we will usually take to be the standard monomials of the $R_n$-ideal $R_nI$ with respect to a term order on $R_n$. In our computational setup, the dependence on the choice of a basis enters via a weight vector that is used to define an elimination term order on~$D_n$. Gauge transformations explicitly describe how passing from one $\CC(x)$-basis of $R_n/R_n I$ to another affects the connection matrices of the resulting system.

\Cref{sec:gaugeth} introduces some theoretical background. \Cref{sec:gaugecomp} presents the algorithms required for the computation of the connection matrices and gauge transforms.

\subsection{Theory}\label{sec:gaugeth}
Let $\{s_1,\ldots,s_m\}$ be a $\CC(x_1,\ldots,x_n)$-basis of $R_n/R_nI$. The $s_j$'s can be chosen to be monomials in the $\partial_i$'s, for instance as the standard monomials of a Gröbner basis of $R_nI$, 
see \Cref{sec:Groebner}. W.l.o.g., we may assume $s_1=1$.
For $f\in \Sol(I)$ a solution to~$I$, denote by $F=(1\bullet f,s_2\bullet f,\ldots,s_m\bullet f)^\top$ the vector of functions formed by applying the operators $s_1,\dots, s_m$ to $f$. Then there exist unique matrices $A_1,\ldots,A_n \in \operatorname{Mat}_{m\times m}(\CC(x_1,\ldots,x_n))$ s.t.\ 
\begin{align}\label{eq:PfaffSys}
\partial_i \bullet F \,=\, A_i \cdot F \, , \quad i\,=\,1,\ldots, n
\end{align}
for \underline{any} $f\in \Sol(I)$. {In \cite{SST00}, the system in \eqref{eq:PfaffSys} is called the {\em Pfaffian system} of $I$ (with respect to the chosen basis). We also refer to the $A_i$'s as the {\em connection matrices} of~$I$.} They encode the transformation on $R_n/R_n I$ given by left-multiplication with $\partial_i$ for the chosen $\CC(x)$-basis. Note, however, that this transformation is \underline{not} $\CC(x)$-linear, but rather extends according to the Leibniz rule.
By construction, the connection matrices fulfill the integrability conditions,~i.e.,
\begin{align}\label{eq:integr}
    \partial_i \bullet A_j-\partial_j \bullet A_i\,=\, \left[A_i,A_j\right] \quad  \text{ for all } \ i,j=1,\ldots,n \, ,
\end{align}
where entry-wise differentiation of the matrices is meant.
Changing basis to $\widetilde{F}=gF$ via some $g\in \GL_m(\CC(x_1,\ldots,x_n))$  yields the system 
$\partial_i \bullet \widetilde{F}= \widetilde{A}_i\cdot \widetilde{F}$, with 
\begin{align}\label{eq : change of var}
\widetilde{A}_i \,=\, gA_ig^{-1}+\left(\partial_i\bullet g\right)g^{-1}, \qquad i=1,\ldots,n.
\end{align} 
This transformation of the connection matrices is called a {\em gauge transformation}, and $\widetilde{A}_i$ is the {\em gauge transform} of $A_i$ (under the {\em gauge matrix} $g$).
Dually---and by using the tensor-hom adjunction---one can equivalently write $D$-ideals in ``connection form'' in terms of a single matrix of differential one-forms which keeps track of all of the $A_i$'s simultaneously. Keeping the same basis, one associates the matrix of differential one-forms $A=A_1\d x_1+\cdots + A_n \d x_n$, the {\em connection matrix} of~$I$. This geometric flavor arises from the fact that, in the holonomic case, a $D_n$-module $D_n/I$ has an underlying vector bundle structure, cf.~\cite{HTT08}, and $\d - A\wedge \, $ defines a flat (also called ``integrable'') connection $\nabla$ on its dual vector bundle (to be precise, at the stalk of the generic point),  
whose flat sections correspond precisely to the solutions of~$I$.

\smallskip
In applications, the considered $D$-ideals often depend on additional parameters, such as a ``small'' parameter~$\varepsilon$---and hence so do the resulting connection matrices.
We say that a connection matrix~$A$ is in {\em $\varepsilon$-factorized form} if ${\varepsilon}^{-1}A $  (or, more generally, ${\varepsilon}^{k}A $ for some integer $k$) is independent of $\varepsilon$.
This form is especially helpful in the context of dimensional regularization of Feynman integrals in particle physics, as this allows for the construction of formal power series solutions in the variable $\varepsilon$ of such systems via the ``path-ordered exponential formalism,'' reducing the computational effort.

\subsection{Computation}\label{sec:gaugecomp}
To compute the connection matrices, we proceed as follows. First, we calculate a Gröbner basis $G$ of the $D_n$-ideal $I$ with respect to an elimination term order $\prec_{(0,v)}$ on $D_n$ with $v \in \mathbb{R}^n_{>0}$ positive. Then, by
\Cref{prop1413}, $G$ is a Gröbner basis of $R_n I$ with respect to the term order $\prec'_{(0,v)}$ on $R_n$. { While the algorithms presented here would work for any elimination order, our presentation and implementation fixes $\prec$ to be the lexicographic order on~$D_n$.}

Later in this section, we will recall an algorithm to reduce elements modulo our Gröbner basis of $R_n I$. Applying it to the  operators $\partial_i s_j$ results in the {\em normal form} of $\partial_i s_j$  w.r.t.\ $G$, where $i=1,\dots,n$ and $j = 1,\dots,m$.
This normal form is of the shape
\begin{align}
    a_{j1}^{(i)}s_1+a_{j2}^{(i)}s_2+\cdots+a_{jm}^{(i)}s_m \, ,
\end{align}
where the coefficients $a_{jk}^{(i)}$ are rational functions in $x_1,\dots,x_n$. Since $G$ is a Gröbner basis with respect to a term order, the normal form is unique, see~\cite[p.~8]{SST00}. Therefore, we can write
\begin{align}
    \partial_i s_j \,=\, \sum_{k = 1}^m a_{jk}^{(i)}s_k + Q^{(i)}_j
\end{align} 
with $Q^{(i)}_j\in R_nI$. Hence $a_{jk}^{(i)}$ is the $(j,k)$-th entry of the matrix $A_i$. 
\smallskip

\Cref{alg:normal form algorithm}, presented below, is an adaptation  of the normal form algorithm given in~\cite[p.~7]{SST00} to the rational Weyl algebra. We can represent an element $P \in R_n$ as
\begin{align*}
    P \,=\, P_\beta \partial^\beta + \text{ lower order terms with respect to } \prec'_{(0,v)} \, ,
\end{align*}
where $P_\beta \in \mathbb{C}(x)$. Similarly, we can represent an element $Q \in D_n$ as 
\begin{align*}
    Q \,=\, Q_{b} \partial^b + \text{lower order terms with respect to } \prec'_{(0,v)} \, ,
\end{align*}
where $Q_{b} \in \mathbb{C}[x]$. We call $\lt_{\prec'_{(0,v)}}(P)\coloneqq P_\beta \partial^\beta$ and $\lt_{\prec'_{(0,v)}}(Q) \coloneqq  Q_{b} \partial^b$ the {\em leading term} of~$P$ and~$Q$, respectively. In contrast to the initial monomials as introduced in \Cref{sec:Groebner}, the leading terms of differential operators are again elements of the (rational) Weyl algebra and, moreover, they contain the coefficients. Note that we regard $Q$ as an element of~$R_n$ and use the restricted order $\prec'_{(0,v)}$ in $R_n$ and not the order $\prec_{(0,v)}$ in $D_n$, even though $Q \in D_n$. The reason for this is the dependence of the leading term on the considered Weyl algebra; it can differ when passing from $D_n$ to $R_n$, as the next example demonstrates.
\begin{example}\label{leading term dependence on space}
    Let $v = (2,1)$ and $\prec_{(0,v)}$ be the elimination 
    term order on~$D_2$ with $\prec$ being the lexicographic order. 
    Let 
    \begin{align*}
        I \,=\, \left\langle x\partial_x^2 - y \partial_y^2+\partial_x -\partial_y , x\partial_x + y\partial_y +1\right\rangle  \eqqcolon \left\langle P_1,P_2\right\rangle .
    \end{align*}
    A Gröbner basis of the $D_2$-ideal $I$ with respect to $\prec_{(0,v)}$ is given by 
    \begin{align*}
        \left\{ y\partial_x\partial_y+\partial_x+y\partial_y^2+\partial_y,\, 
    x\partial_x+y\partial_y+1, \, 
    xy\partial_y^2-y^2\partial_y^2+x\partial_y-3y\partial_y-1\right\} \, .
    \end{align*}
    By \Cref{prop1413}, it is also a Gröbner basis of $R_2 I$ with respect to~$\prec'_{(0,v)}$. The third Gröbner basis element, $G_3$, as an element of~$D_2$ has the leading term $\operatorname{lt}_{\prec_{(0,v)}}(G_3) = xy\partial_y^2$, but considered as an element of $R_2$, its leading term is $\operatorname{lt}_{\prec'_{(0,v)}}(G_3) = (x-y)y\partial_y^2$.
\end{example}
We will return to this example in \Cref{sec:implementations} to demonstrate the methods implemented in our package.

\smallskip

If, in the notation as above, we have $\beta_i \geq b_i$ for all $1 \leq i \leq n$, i.e., if the initial monomial of~$Q$ divides the initial monomial of $P$ in~$R_n$, the reduction of $P$ by $Q$ is defined by
\begin{align}\label{eq:reduction} 
    \operatorname{red}_{\prec'_{(0,v)}}(P,Q) \, \coloneqq \, P-\frac{P_\beta}{Q_{b}}\, \partial^{\beta-b}Q \, . 
\end{align}
In this case, it coincides with the S-pair of $P$ and $Q$, see \cite[p.~7]{SST00}. Observe that we multiplied $Q$ by $(P_\beta / Q_{b})\partial^{\beta-b}$ to cancel the leading terms of $P$ and~$Q$. It leads to the following algorithm.

\begin{algorithm}[H]
    \caption{(Normal form algorithm in the rational Weyl algebra){\bf .}}\label{alg:normal form algorithm}
    \begin{algorithmic}
    \smallskip
        \Require $P \in R_n$, $v \in \mathbb{R}^n_{>0}$, a Gröbner basis $G$ of a $D_n$-ideal $I$ with respect to $\prec_{(0,v)}$ on $D_n$, 
        \linebreak \hspace*{9.4mm} with $\prec$ being the lexicographic order on $D_n$. 
        \Ensure The {\tt normalForm} of $P$ by $G$ in $R_n$ with respect to $\prec'_{(0,v)}$.
        \smallskip
        \If{$P==0$}
            \State \Return $P$
        \EndIf
        \While{$\exists G_i \in G$ s.t.\ $\init_{\prec'_{(0,v)}}(G_i) \mid \init_{\prec'_{(0,v)}}(P)$}
            \State $P \coloneqq  \operatorname{red}_{\prec'_{(0,v)}}(P,G_i)$
        \EndWhile \smallskip
        \State \Return lt$_{\prec'_{(0,v)}}(P)$ $+$ \texttt{normalForm}($P-\leadTerm_{\prec'_{(0,v)}}(P)$, $G$)
    \end{algorithmic}
\end{algorithm}

Note that \Cref{alg:normal form algorithm} terminates since the leading term becomes successively smaller and $1$ is the smallest monomial of any term order $\prec$ on $D_n$. Note also that \Cref{prop1413} ensures that the normal form of any $P \in R_n I$ is $0$. We point out that, in general, it is not sufficient to carry out the reduction in the {\tt while} loop only once for the elements of the Gröbner basis.

\smallskip
As seen above, the normal form algorithm allows us to compute the connection matrices of an ideal $I \subset D_n$, which we summarize in pseudo-code in the next algorithm.

\begin{algorithm}[H]
    \caption{(Connection matrices with respect to standard monomials){\bf .}} \label{alg:connectionMatrices}
    \begin{algorithmic} 
    \smallskip
        \Require A $D_n$-ideal $I$ of finite holonomic rank $m <\infty$ and a positive vector $v \in \RR_{>0}^n$.
        \Ensure The {\tt connectionMatrices} $A_1,\dots,A_n \in \CC(x)^{m \times m}$ of $I$ with respect to the \linebreak \hspace*{13.4mm} standard monomials for~$\prec'_{(0,v)}$.\smallskip 
        \State $G\coloneqq$  Gröbner basis of the $D_n$-ideal $I$ with respect to $\prec_{(0,v)}$
        \State $\{s_1 \prec'_{(0,v)} s_2 \prec'_{(0,v)} \dots \prec'_{(0,v)} s_m\}\coloneqq \{\partial^\beta \colon \beta \in \NN^n \text{ s.t.\ } \init_{\prec'_{(0,v)}}(P) \nmid \xi^\beta \ \forall P \in G\}$ 
        \For{$i$ from $1$ to $n$}
          \For{$j$ from $1$ to $m$}
            \State $P \coloneqq $ \texttt{normalForm}($\partial_i s_j$,\, $G$) \textcolor{black!60}{\Comment{normal form computation in $R_n$ w.r.t.\ $\prec'_{(0,v)}$}}
            \For{$k$ from $1$ to $m$}
               \State $a_{jk}^{(i)} \coloneqq  $ coefficient of the monomial $s_k$ in $P$
            \EndFor
          \EndFor
          \State $A_i \coloneqq (a_{jk}^{(i)}) \in \CC(x)^{m \times m}$
        \EndFor
        \State \Return $A_1,\dots,A_n$.
    \end{algorithmic}
\end{algorithm}
\vspace*{-2mm}

The following example shows how to compute the connection matrices via this algorithm.
\begin{example} 
    Let $I$ be the $D_2$-ideal from \Cref{leading term dependence on space} and again $v=(2,1)$. We already determined a Gröbner basis $G$ of $I$ and $R_2I$ with respect to $\prec_{(0,v)}$ and $\prec'_{(0,v)}$, respectively. We have $\init_{\prec'_{(0,v)}}(R_2 I) = \langle \xi_x \xi_y , \xi_x, \xi_y^2 \rangle$. The holonomic rank of $I$ is the number of monomials~$\xi^\alpha$ that are not contained in the initial ideal $\init_{\prec'_{(0,v)}}(R_2 I)$.
    Therefore, $\rank(I) = 2$, with $s_1=1$ and $s_2=\partial_y$ being the standard monomials. We hence choose $\{1,\partial_y\}$ as our \mbox{$\mathbb{C}(x,y)$-basis} of~$R_2/R_2 I$. Following the above algorithm, we have $\operatorname{red}_{\prec'_{(0,v)}}(\partial_x, x\partial_x+y\partial_y+1) = \partial_x-\frac{1}{x}(x\partial_x+y\partial_y+1) = -\frac{y}{x}\partial_y-\frac{1}{x}$. Consequently, we get the normal form
    \begin{align*}
        \texttt{normalForm}(\partial_x s_1, G) &\,=\, \texttt{normalForm}(\partial_x, G) \,=\, -\frac{y}{x}\partial_y-\frac{1}{x} \,.
    \end{align*}
    Analogously, we obtain
    \begin{align*}
        \texttt{normalForm}(\partial_x s_2, G) &\,=\, \texttt{normalForm}(\partial_x\partial_y, G) \,=\, -\frac{x+y}{x(x-y)} \partial_y - \frac{1}{x(x-y)} \, , \\
        \texttt{normalForm}(\partial_y s_1, G) &\,=\, \texttt{normalForm}(\partial_y, G) \,=\, \partial_y \, ,\\
        \texttt{normalForm}(\partial_y s_2, G) &\,=\, \texttt{normalForm}(\partial_y^2, G) \,=\, \frac{3y-x}{(x-y)y}\partial_y + \frac{1}{(x-y)y} \, .
    \end{align*}
    From the first and last two normal forms, we obtain the connection matrices of $I$ as
    \renewcommand\arraystretch{1.2}
    \begin{align*}
        A_1 = \begin{pmatrix}
            -\frac{1}{x} & - \frac{y}{x}\\
            - \frac{1}{x(x-y)} & -\frac{x+y}{x(x-y)}
        \end{pmatrix} \quad \text{and} 
        \quad A_2 = \begin{pmatrix}
            0 & 1\\
            \frac{1}{(x-y)y} & \frac{3y-x}{(x-y)y}
        \end{pmatrix},
    \end{align*}
and the corresponding connection matrix as 
    \begin{align*}
        \qquad \qquad
        A \,=\,  A_1 \d x + A_2 \d y \,=\, \begin{pmatrix}
            -\frac{1}{x} \d x & - \frac{y}{x} \d x + \d y\\
            - \frac{1}{x(x-y)} \d x + \frac{1}{(x-y)y} \d y & -\frac{x+y}{x(x-y)}\d x + \frac{3y-x}{(x-y)y} \d y
        \end{pmatrix}. \qquad \qquad \qedhere
    \end{align*}
    \renewcommand\arraystretch{1}
\end{example}

We can now compute the connection matrices in arbitrary $\CC(x)$-bases of $R_n/R_nI$ (not necessarily given by standard monomials) by first computing the connection matrices~$A_i$ using \Cref{alg:connectionMatrices} and then computing the resulting gauge transforms~$\widetilde{A}_i$. 

\begin{algorithm}[H]
    \caption{(Gauge transformation to an arbitrary basis){\bf .}}\label{alg:gaugeTransform}
    \begin{algorithmic} 
        \Require A $D_n$-ideal $I$ of $\rank(I)=m < \infty$, a $\CC(x)$-basis $\{r_1,\dots,r_m\}$ of $R_n/R_n I$,  $v \in \RR_{>0}^n$.
        \Ensure The connection matrices of $I$ in the basis $\{r_1,\dots,r_m\}$. \smallskip
        \State $A\coloneqq {\tt connectionMatrices}(I, v)$ 
        \State $G\coloneqq$  Gröbner basis of the $D_n$-ideal $I$ with respect to $\prec_{(0,v)}$
        \State $\{s_1 \prec'_{(0,v)} s_2 \prec'_{(0,v)} \dots \prec'_{(0,v)} s_m\} \,\coloneqq\, \{\partial^\beta \colon \beta \in \NN^n \text{ s.t.\ } \init_{\prec'_{(0,v)}}(g) \nmid \xi^\beta \ \forall g \in G\}$ 
        \For{$j$ from $1$ to $m$}
          \State $P \coloneqq $ \texttt{normalForm}($r_j$,\, $v$,\, $G$) \textcolor{black!60}{\Comment{normal form computation in $R_n$ w.r.t.\ $\prec'_{(0,v)}$}}
          \For{$k$ from $1$ to $m$}
            \State $g_{jk} \coloneqq $ coefficient of the monomial $s_k$ in $P$.
          \EndFor
        \EndFor
        \State $g \coloneqq (g_{jk}) \in \CC(x)^{m \times m}$. 
        \smallskip
        \State \Return $\widetilde{A}_i \coloneqq gA_ig^{-1}+(\partial_i\bullet g)g^{-1}$ for  $i=1,\ldots,n$.
    \end{algorithmic}
\end{algorithm}
\vspace*{-1mm}
\noindent Examples for gauge transformations are provided in the documentation of our package.

\section{Implementation}\label{sec:implementations}
The core methods implemented in our package are the following: {\tt normalForm},\linebreak {\tt standardMonomials}, {\tt connection{Form}} {(with {\tt connectionMatrix} as an alias),\linebreak {\tt pfaffianSystem} (with {\tt connectionMatrices} as an alias)}, {\tt isIntegrable}, {\tt gaugeMatrix}, {\tt gaugeTransform}, and {\tt isEpsilonFactorized}.
In addition, we extended the definition of {\tt makeWeylAlgebra} to be able to set weight orders that are refined by the lexicographic order built upon $\partial_1 \succ \cdots \succ \partial_n \succ x_1 \succ \cdots \succ x_n$.

All of our described algorithms are fundamentally based on the computation of {\tt normalForm} in the rational Weyl algebra $R_n$. However, at the time of submission of the present article, rational Weyl algebras are not yet implemented in {\em Macaulay2}.\linebreak As a stand-in, we use $\CC(x)[\xi]$ to internally represent operators from $R_n$ in their standard form. In particular, the implementation of the reduction step~\eqref{eq:reduction} computes $\partial^{\beta-b}Q$ in $D_n$ and substitutes it into $\CC(x)[\xi]$ only afterwards.
Additionally, the connection matrices of ideals in a Weyl algebra {\tt D} should be defined over the fraction field of the underlying polynomial ring. We refer to this as {\tt baseFractionField(D)}. 

We demonstrate the usage of all of these commands on-the-fly for the \mbox{$D_2$-ideal}\linebreak $I=\langle x\partial_x^2-y\partial_y^2+\partial_x-\partial_y, \, x\partial_x+y\partial_y+1 \rangle$ from \Cref{leading term dependence on space}.
We found the output of this simple example to be well-suited for displaying purposes. More interesting examples can be found on the MathRepo page of our project. There, we revisit an annihilating $D$-ideal of a correlation function in cosmology as in \cite[(11)]{FPSW24} as well as an annihilating $D$-ideal of a massless one-loop triangle Feynman integral~\cite{HPSZ} from particle physics.

\subsection{Connection matrices}
The following functions allow us to define Weyl algebras endowed with the term order $\prec_{(0,v)}^{\text{lex}}$ with positive $v\in \RR_{>0}^n$ and to compute connection matrices of $D$-ideals in the basis of the standard monomials of a Gröbner basis of $R_nI$ with respect to the term order~${\prec_{(0,v)}^{\text{lex},'}} $ on~$R_n$.

\vspace*{-3mm}
\subsubsection {\tt makeWeylAlgebra(R,v)} defines the Weyl algebra of a polynomial ring {\tt R}. The obtained Weyl algebra is endowed with the term order $\prec_{(0,{\tt v})}^{\text{lex}}$, with (0,{\tt v}) being an element of $ 0\times \RR^n_{>0}\subset \mathcal{W}$, assigning weight $0$ to the $x_i$'s and positive weight $v_i$ to $\partial_i$. This guarantees that the resulting term order on~$D_n$ is an elimination term order. 
\vspace*{1mm}
\begin{small}
\begin{lstlisting}[language=Macaulay2, frame=single]
i1 : needsPackage "ConnectionMatrices";
i2 : D = makeWeylAlgebra(QQ[x,y], {2,1});
\end{lstlisting}
\end{small}

\medskip
\noindent {\tt makeWeylAlgebra(R)} defines the Weyl algebra of {\tt R} endowed with the term order~$\prec^{\text{lex}}$. 
This is equivalent to setting {\tt v = \{0,...,0\}} above.
Note that the polynomial ring used in the definition of the Weyl Algebra is allowed to have coefficients in the fraction field of a polynomial ring, such as $\CC(\varepsilon)$. This is often necessary when considering Weyl algebras arising from physics. An example of these two cases is as follows.

\begin{small}
\vspace*{1mm}
\begin{lstlisting}[language=Macaulay2, frame=single]
i3 : Reps = frac(QQ[eps])[x];
i4 : Deps = makeWeylAlgebra(Reps);
\end{lstlisting}
\end{small}

\vspace*{-2mm}
\subsubsection{\tt normalForm(P,Q)} computes the normal form of {\tt P} with respect to {\tt Q}. Both ${\tt P}$ and ${\tt Q}$ have to be elements of~$D_n$. The reduction step is carried out in the rational Weyl algebra.
\begin{small}
\vspace*{1mm}
\begin{lstlisting}[language=Macaulay2, frame=single]
i5 : use D;
i6 : P = dx;
i7 : Q = x*dx + y*dy + 1;
i8 : normalForm(P,Q)
     -y      -1
o8 = --*dy + --
      x       x
o8 : frac(QQ[x..y])[dx, dy]
\end{lstlisting}
\end{small}

\vspace*{-2mm}
\subsubsection*{\tt normalForm(P,G)} computes the normal form of {\tt P}, an element of the Weyl algebra $D_n$, with respect to a list {\tt G} of elements in the Weyl algebra (typically a Gröbner basis of a $D_n$-ideal). 
\begin{small}
\vspace*{1mm}
\begin{lstlisting}[language=Macaulay2, frame=single]
i9 : I = ideal(x*dx^2-y*dy^2+dx-dy,x*dx+y*dy+1)
                2            2
o9 = ideal (x*dx  + dx - y*dy  - dy, x*dx + y*dy + 1)
o9 : Ideal of D
i10 : G = gens gb I
o10 = | xydy^2-y2dy^2+xdy-3ydy-1 xdx+ydy+1 ydxdy+dx+ydy^2+dy |
i11 : normalForm(dx*dy, flatten entries G)
       - x - y         -1
o11 = --------*dy + --------
       2             2
      x  - x*y      x  - x*y
o11 : frac(QQ[x..y])[dx, dy]
\end{lstlisting}
\end{small}

\vspace*{-2mm}
\subsubsection{\tt standardMonomials(I)} computes the standard monomials of the Gr\"obner basis of~$R_n${\tt I} with respect to the restriction of the chosen term order on~$D_n$ to~$R_n$.
\vspace*{1mm}
\begin{small}
\begin{lstlisting}[language=Macaulay2, frame=single]m
i12 : m = holonomicRank I
o12 = 2
i13 : SM = standardMonomials I
o13 = {1, dy}
\end{lstlisting}
\end{small}

\medskip
    In order to compute the standard monomials with respect to a term order for another weight vector, it is necessary to define a new Weyl algebra with that term order.
    \vspace*{1mm}
\begin{small}
\begin{lstlisting}[language=Macaulay2, frame=single]
i14 : D2 = makeWeylAlgebra(QQ[x,y],{1,2});
i15 : SM2 = standardMonomials(sub(I,D2))
o15 = {1, dx}
\end{lstlisting}
\end{small}

\vspace*{-2mm}
\subsubsection{\tt pfaffianSystem(I)} 
computes the list of connection matrices {\tt A} of the $D_n$-ideal~{\tt I} 
w.r.t.\ the standard monomials of a Gröbner basis of $R_nI$ for the chosen term order on~$D_n$.
\vspace*{1mm}
\begin{small}
\begin{lstlisting}[language=Macaulay2, frame=single]
i16 : A = pfaffianSystem sub(I,D)
o16 = {| (-1)/x       (-y)/x         |, | 0         1               |}
       | (-1)/(x2-xy) (-x-y)/(x2-xy) |  | 1/(xy-y2) (-x+3y)/(xy-y2) |
\end{lstlisting}
\end{small}

\medskip
\noindent{\tt pfaffianSystem(I,B)} computes the connection matrices of the $D_n$-ideal {\tt I} w.r.t.\ a chosen $\CC(x)$-basis {\tt B} of $R_n/R_nI$; below, we chose {\tt B} to be the standard monomials {\tt SM2} from~{\tt i16}. 
\vspace*{1mm}
\begin{small}
\begin{lstlisting}[language=Macaulay2, frame=single]
i17 : A2 = pfaffianSystem(I,SM2)
o17 = {| 0            1               |, | (-1)/y    (-x)/y        |}
       | (-1)/(x2-xy) (-3x+y)/(x2-xy) |  | 1/(xy-y2) (x+y)/(xy-y2) |
\end{lstlisting}
\end{small}

\vspace*{-2mm}
\subsubsection{\tt connectionForm(I)} displays the connection matrix of the $D_n$-ideal {\tt I}. \smallskip\\{\bf Nota bene:} This command is to be used for displaying purposes only; this matrix is \underline{not} encoded in the respective ring of differential forms and should therefore not be used for further computations in {\em Macaulay2}.
\vspace*{1mm}
\begin{small}
\begin{lstlisting}[language=Macaulay2, frame=single]
i18 : connectionForm(I)
o18 = | (-1)/xdx                   (-y)/xdx+dy                        |
      | (-1)/(x2-xy)dx+1/(xy-y2)dy (-x-y)/(x2-xy)dx+(-x+3y)/(xy-y2)dy |
\end{lstlisting}
\end{small}

\vspace*{-2mm}
\subsubsection{\tt isIntegrable(\{A\_1,...,A\_n\})} checks whether a given list {\tt \{A\_1,...,A\_n\}} of matrices  over (the fraction field of) a polynomial ring fulfills the integrability conditions~\eqref{eq:integr}. 
\begin{small}
\begin{lstlisting}[language = Macaulay2,frame=single]
i19 : isIntegrable(A)
o19 = true
\end{lstlisting}
\end{small}

\subsection{Gauge transformation}
Given the connection matrices of a $D_n$-ideal $I$ with respect to a given basis, it is possible to rewrite them with respect to another basis via the gauge transform~\eqref{eq : change of var}. This is implemented in our package via the following commands. 

\vspace*{-2mm}
\subsubsection{\tt gaugeMatrix(I,B)} outputs the matrix that encodes the gauge matrix from the basis consisting of the standard monomials of the Gr\"obner basis of the $R_n$-ideal generated by {\tt I} to the basis {\tt B} (here {\tt \{1,dx\}}).
\begin{small}
\vspace*{1mm}
\begin{lstlisting}[language=Macaulay2, frame=single]
i20 : F = baseFractionField(D);
i21 : g = gaugeMatrix(I,{1_D,dx_D})
o21 = | 1      0      |
      | (-1)/x (-y)/x |
               2      2
o21 : Matrix F  <-- F
\end{lstlisting}
\end{small}

\vspace*{-2mm}
\subsubsection{\tt gaugeTransform(g,A)} computes the gauge transform of the list of connection matrices {\tt A} of a $D$-ideal {\tt I} for the gauge matrix {\tt g} via~\eqref{eq : change of var}. This results in the list {\tt A2'} of gauge-transformed matrices. For example: 
\vspace*{1mm}
\begin{small}
\begin{lstlisting}[language=Macaulay2, frame=single]
i22 : A2' = gaugeTransform(g,A)
o22 = {| 0            1               |, | (-1)/y    (-x)/y        |}
       | (-1)/(x2-xy) (-3x+y)/(x2-xy) |  | 1/(xy-y2) (x+y)/(xy-y2) |
i23 : A2 == A2'
o23 = true
\end{lstlisting}
\end{small}

\vspace*{-2mm}
\subsubsection{\tt isEpsilonFactorized(A,eps)} checks whether a family of connection matrices {\tt A} is in {\tt eps}-factorized form, that is, if it is possible to factor out a power of the variable {\tt eps} so that the remaining matrix is independent of {\tt eps}. 
\vspace*{1mm}
\begin{small}
\begin{lstlisting}[language=Macaulay2, frame=single]
i24 : use Deps;
i25 : I = ideal(x*(1-x)*dx^2 - eps*(1-x)*dx);
i26 : B = {sub(1,Deps),sub(1/eps,Deps)*dx};
i27 : Aeps = connectionMatrices(I,B)
o27 = {| 0 eps   |}
       | 0 eps/x |
i28 : isEpsilonFactorized(Aeps,eps)
o28 = true
\end{lstlisting}
\end{small}

\bigskip \bigskip
\noindent {\bf Acknowledgments.}
This project began at the workshop {\it Macaulay2 in the Sciences} held at MPI-MiS Leipzig in November 2024. We are grateful to Devlin Mallory and Carlos Gustavo Rodriguez Fernandez, who contributed to the development of this package during the workshop, and to Anton Leykin, Mike Stillman, and Bernd Sturmfels for helpful discussions.

\bigskip
\noindent {\bf Funding statement.}
The research of ALS and NW is funded by the European Union (ERC, UNIVERSE PLUS, 101118787). Views and opinions expressed are, however, those of the author(s) only and do not necessarily reflect those of the European Union or the European Research Council Executive Agency. Neither the European Union nor the granting authority can be held responsible for them.

\addcontentsline{toc}{section}{References}
\bibliographystyle{abbrv} 

\end{document}